\documentclass[journal]{IEEEtran}
\usepackage{ifpdf}

\usepackage{algorithm}
\usepackage{algpseudocode}
\usepackage{color}
\usepackage{url}
\ifpdf
\usepackage[pdftex]{graphicx}
\else
\usepackage[dvips]{graphicx}
\fi
\usepackage{tabularx}
\DeclareGraphicsExtensions{.eps}
\usepackage[caption=false,font=footnotesize]{subfig}
\usepackage{amsmath,amsthm,amsfonts,amssymb}
\usepackage[noadjust]{cite}
\usepackage{nomencl}
\usepackage{enumerate}

\newcommand{\ud}[1]{_\mathrm{#1}}
\newcommand{\up}[1]{^\mathrm{#1}}
\newcommand\Tstrut{\rule{0pt}{2.6ex}}         
\newcommand\Bstrut{\rule[-0.9ex]{0pt}{0pt}}   

\makeatletter
\newcommand{\subalign}[1]{%
	\vcenter{%
		\Let@ \restore@math@cr \default@tag
		\baselineskip\fontdimen10 \scriptfont\tw@
		\advance\baselineskip\fontdimen12 \scriptfont\tw@
		\lineskip\thr@@\fontdimen8 \scriptfont\thr@@
		\lineskiplimit\lineskip
		\ialign{\hfil$\m@th\scriptstyle##$&$\m@th\scriptstyle{}##$\crcr
			#1\crcr
		}%
	}
}
\makeatother

\newcommand{\sumtb}{\sum_{\subalign{t&\in T \\ b&\in B}}}
\newcommand{\sumti}{\sum_{\subalign{t&\in T \\ i&\in I}}}

\begin{document}

\title{Scalable Planning for Energy Storage \\in  Energy and Reserve Markets }	

\author{Bolun~Xu,~\IEEEmembership{Student Member,~IEEE,}
	Yishen~Wang,~\IEEEmembership{Student Member,~IEEE,}
	Yury~Dvorkin,~\IEEEmembership{Member,~IEEE,}
	Ricardo~Fern\'andez-Blanco,
	C. A. Silva-Monroy, \IEEEmembership{Member,~IEEE,}
	Jean-Paul Watson, \IEEEmembership{Member,~IEEE,}\\
	and Daniel~S.~Kirschen,~\IEEEmembership{Fellow,~IEEE}
	\vspace{-.7cm}
	\thanks{This work was supported by the US Department of Energy under grant  1578574.}
	\thanks{B.~Xu, Y.~ Wang, R.~Fern\'andez-Blanco, and D.~S.~Kirschen are with University of Washington, USA (emails: \{xubolun, ywang11, kirschen\}@uw.edu and ricardo.fcarramolino@gmail.com). }
	\thanks{Y.~Dvorkin is with New York University, USA (email: dvorkin@nyu.edu). }
	\thanks{C.~Silva-Monroy and J.~Watson are with Sandia National Laboratories (emails: \{casilv, jwaston\}@sandia.gov). }%
}  

\maketitle
\makenomenclature

\begin{abstract}

Energy storage can facilitate the integration of renewable energy resources by providing arbitrage and ancillary services. Jointly optimizing energy and ancillary services in a centralized electricity market reduces the system's operating cost and enhances the profitability of energy storage systems. However, achieving these objectives requires that storage be located and sized properly. We use a bi-level formulation to optimize the location and size of energy storage systems which perform energy arbitrage and provide regulation services. Our model also ensures the profitability of investments in energy storage by enforcing a rate of return constraint. Computational tractability is achieved through the implementation of a primal decomposition and a subgradient-based cutting-plane method. We test the proposed approach on a 240-bus model of the Western Electricity Coordinating Council (WECC) system and analyze  the effects of different storage technologies, rate of return requirements, and regulation market policies on ES participation on the optimal storage investment decisions.
We also demonstrate that the proposed approach outperforms exact methods in terms of solution quality and computational performance.

\end{abstract}


\IEEEpeerreviewmaketitle

\section*{Nomenclature}
\addcontentsline{toc}{section}{Nomenclature}

\subsection{Sets and Indices}

\begin{IEEEdescription}[\IEEEusemathlabelsep\IEEEsetlabelwidth{$V_1,V_2,V_3$}]
\item[$B, B\up{E}, B\up{N}$] Set of buses and subset of buses with and without energy storage, indexed by $b$
\item[$I, I_b$] Set of generators and subset of generators connected to bus $b$, indexed by $i$
\item[$J$] Set of typical days, indexed by $j$
\item[$L$] Set of transmission lines, indexed by $l$
\item[$T$] Set of time intervals, indexed by $t$
\item[$o(l), r(l)$] Sending and receiving ends of line $l$
\end{IEEEdescription}

\vspace{-5mm}
\subsection{Variables}
\begin{IEEEdescription}[\IEEEusemathlabelsep\IEEEsetlabelwidth{$V_1,V_2,V_3$}]
\item[$C\up{[\cdot]}(\cdot)$] Cost function for different problems, \$
\item[$e\up{R}_b, p\up{R}_{b}$] Energy storage energy (MWh) and power (MW) rating.
\item[$e\up{soc}_{j,t,b}$] Energy storage state of charge, MWh
\item[$f_{j,t,l}$] Line power flow, MW
\item[$p\up{ch}_{j,t,b}, p\up{dis}_{j,t,b}$] Energy storage charging and discharging rates, MW
\item[$p\up{g}_{j,t,i}$] Generator output power, MW
\item[$p\up{rs}_{j,t,b}$] Renewable spillage, MW
\item[$r\up{ed}_{j,t,b}, r\up{eu}_{j,t,b}$] Regulation down and up provided by energy storage, MW
\item[$r\up{gd}_{j,t,i},r\up{gu}_{j,t,i}$] Regulation down and up provided by generator, MW
\item[$\mathnormal{x}\up{[\cdot]}$] Vector of decision variables
\item[$y$] Cutting-plane method decision variable vector
\item[$\theta_{j,t,b}$] Voltage phase angle
\item[$\lambda\up{lmp}_{j,t,b}$] Locational marginal price, \$/MWh
\item[$\lambda\up{rd}_{j,t}, \lambda\up{ru}_{j,t}$] Price for down/up regulation, \$/MWh
\item[$\varphi\up{[.]}_{[.]}, \psi\up{[.]}_{[.]}, \gamma\up{[.]}_{[.]}$] Dual variables associated with upper bounds, lower bounds, and equality constraints
\end{IEEEdescription}

\subsection{Parameters}

\begin{IEEEdescription}[\IEEEusemathlabelsep\IEEEsetlabelwidth{$V_1,V_2,V_3$}]
\item[$c\up{p}, c\up{e}$] Daily prorated capital costs of energy storage, \$/MW and \$/MWh
\item[$c\up{g}_i,c\up{gd}_i,c\up{gu}_i$] Hourly incremental, regulation up, and regulation down costs of generators, \$/MWh
\item[$c\up{dis},c\up{ch}$] Hourly incremental discharge and charge cost of ES, \$/MWh
\item[$c\up{eu},c\up{ed}$] Hourly regulation up and down cost of ES, \$/MWh
\item[$c\up{rs}$] Value of renewable spillage, \$/MWh
\item[$D_{j,t,b}$] Load demands, MW
\item[$F_l$] Capacity of transmission lines, MW
\item[$G\up{max}_i,G\up{min}_i$] Maximum and minimum power production of generators, MW
\item[$G\up{rn}_{j,t,b}, G\up{rs}_{j,t,b}$] Renewable power maximum expected forecast and maximum allowable spillage, MW
\item[$R\up{u}_i, R\up{d}_i$] Ramp down and ramp up capacity of generators, MW/h
\item[$T\up{ru}, T\up{rd}$] Ramp down/up speed requirement for down/up regulation, h
\item[$T\up{es}$] Continuous full dispatch time requirement for energy storage, h
\item[$X_l$] Transmission line reactance

\item[$\omega_j$] Weight of typical day
\item[$\rho\up{min}, \rho\up{max}$] Minimum and maximum allowable values of the power/energy ratio of storage systems, h$^{-1}$
\item[$\eta\up{ch}, \eta\up{dis}$] Charging and discharging efficiency of energy storage
\item[$\phi\up{D}, \phi\up{R}$] Regulation requirement as a percentage of the load and renewable injections
\item[$\nu$] Iteration index
\item[$\epsilon$] Relative tolerance for system cost savings
\end{IEEEdescription}

\section{Introduction}

Energy storage (ES) is a highly flexible resource that has the potential to facilitate the integration of renewable energy sources such as wind and solar~\cite{barton2004energy,oudalov2006value}. U.S. system operators and regulators have recognized ES as the key technology in achieving sustainability in the power sector~\cite{ercot_es}. For instance, 263~MW of the ES capacity has already been deployed in PJM~\cite{pjm_es}. In ISO New England, 94~MW of the battery energy storage capacity has been proposed for deployment as of January 2016~\cite{isone_outlook}. The California Public Utilities Commission has mandated a merchant ES procurement goal of 1325~MW by 2020~\cite{cal_rm}.

ES can relieve congestion in the transmission system by performing spatio-temporal arbitrage, make possible a more optimal dispatch of conventional generators and hence reduce the cost of dealing with the intermittency of renewable resources~\cite{pandzic2015near}. Advanced ES technologies such as batteries also provide rapid responses in reserve and regulation services, which lower the ancillary service procurement requirements and reduce the cost of handling the stochasticity of wind and solar~\cite{xu2016comparison}. 
However, because it is likely that many of these ES systems will be deployed by private investors, we should not consider only whether they provide a social benefit in terms of reduced operating cost, but also whether they generate a sufficient return on investment~\cite{dvorkin2016ensuring,wang2017look}. Since ES operators may wish to participate in multiple electricity markets to increase their profits~\cite{krishnan2015optimal, xu2014bess, sandia_es, ny_es,akhavan2014optimal}, we need to consider the arbitrage jointly with the provision of ancillary services when identifying opportunities for investments in ES.

An accurate long-term planning decision must account for its impact on short-term system operations~\cite{pudjianto2014whole}. However, solving a single optimization problem that includes the entire planing horizon (i.e., a full-year operation) is far beyond what is computationally tractable at this point in time.
To overcome such computation barriers, heuristic ES planning models \cite{pandzic2015near, dvijotham2014storage} split ES siting and sizing into sequential decisions according to heuristic rules. While heuristic models are solvable over longer planning horizons, they may produce suboptimal planning decisions. To obtain more rigorous planning decisions, stochastic programming has been extensively incorporated in power system planning problems~\cite{gorenstin1993power, buygi2003transmission, de2008transmission,garces2009bilevel}.
Stochastic planning models co-optimize siting and sizing decisions on ES over a set of selected representative scenarios~\cite{buygi2003transmission, de2008transmission, oh2011optimal, bayram2015stochastic, qiu2016stochastic, wogrin2015optimizing,fernandez2016optimal}. The computational complexity of a stochastic planning model depends on the number of scenarios, and a sufficient number of scenarios must be considered for effective representations of the uncertain renewable generation resources and demand. Although a larger number of scenarios improves the robustness of the planning result, such formulated problems can be computationally intractable when applied to large power systems~\cite{winston2004operations}.
Besides, adding additional planning criteria, such as a guarantee of the ES investment payback rates~\cite{dvorkin2016ensuring,nasrolahpour2016strategic}, and a co-optimization of energy and reserve markets~\cite{krishnan2015optimal}, will further increase the complexity of the planning model.

The aforementioned ES planning approaches aim to trade-off modeling accuracy and computational complexity. Still, solving ES planning problems in realistically large systems is a non-trivial task and the modeling accuracy has been sacrificed for the sake of solvability. Nasrolahpour~\emph{et al.}~\cite{nasrolahpour2016strategic} formulated strategic ES sizing in energy markets as a bi-level problem, and adopted a solution algorithm which combines mathematical programming with equilibrium constraints (MPEC) with Benders decomposition. However, their algorithm takes hours to solve the bi-level ES sizing problem on a single bus case study. We propose a decomposition algorithm that provides more accurate and faster solutions to the ES planning problem for a large number of scenarios.
This paper makes the following contributions:
\begin{itemize}
	\item It formulates the optimal ES profit-constrained siting and sizing problem in a joint energy and reserve market as a bi-level problem considering the perspectives of the system operator in anticipation that energy storage would act as profit-seeking entities in a market environment
	\item It describes and tests a solution method which combines primal decomposition with subgradient cutting-planes. This solution method is scalable to any planning scenarios, and has non-heuristic terminating criteria.
	\item It benchmarks the computational performance of this algorithm against an exact linear programming (LP) approach, and demonstrates the accuracy and scalability of this algorithm.
	\item It uses compressed air energy storage and lithium-ion batteries to represent two different types of ES technologies, and compares their investment for different regulation market policies.
	\item It analyzes the effect of a minimum profit constraint on the ES siting and sizing decisions as well as on the system operating cost.
\end{itemize}
All simulations are performed on a modified version of the 240-bus system of Western Electricity Coordinating Council (WECC)~\cite{price2011reduced}. The WECC system is a realistically large testbed for planning studies that demonstrates scalability of the proposed solution method. 
Section~\ref{Sec:PD} formulates the bi-level optimal ES siting and sizing model. Section~\ref{Sec:Meth} proposes the solution method to the formulated problem. Section~\ref{Sec:Case} describes the test-bed system parameters. Section~\ref{Sec:Test} presents and discusses the numerical results. Finally, conclusions are drawn in Section~\ref{Sec:Con}.

\section{Problem Formulation}~\label{Sec:PD}

We formulate the optimal ES siting and sizing as a bi-level problem. The upper-level (UL) problem identifies the ES investment decisions which minimize the overall system cost over a set of typical (or representative) days, while the lower-level (LL) problems minimize the operating cost of each typical day.

\subsection{Upper-Level Problem: Energy Storage Siting and Sizing} 

The UL problem minimizes the total system cost ($C\up{S}$) over all typical days, i.e. the sum of the expected system operating cost and of the ES investment cost:

\begin{gather}\label{PD:UL_Obj}
\textstyle\min_{\mathnormal{x}\up{U}} C\up{S}(x\up{P}_j) :=\sum_{j\in J}\omega_j C\up{P}_j(\mathnormal{x}\up{U}, \mathnormal{x}\up{p}_j)+C\up{E}(\mathnormal{x}\up{U})\,,
\end{gather}
where $\mathnormal{x}\up{U}$ are upper-level decision variables, and $\mathnormal{x}\up{p}_j$ are lower-level decision variables. The system operating cost is the sum of the dispatch cost $C\up{P}_j$ for each typical days weighted by the relative importance $\omega_j$ of the days it represents. The ES investment cost $C\up{E}$ is calculated based on both the power rating $p\up{R}_b$ and the energy capacity $e\up{R}_b$ of the ES installed at each bus $b\in B$:

\begin{gather}
\textstyle C\up{E}(\mathnormal{x}\up{U}) := \sum_{b\in B} (c\up{p} p\up{R}_b + c\up{e} e\up{R}_b)\,,
\end{gather} 
where $\mathnormal{x}\up{U} = \{ p\up{R}_b$, $e\up{R}_b \}$.
This problem is constrained by limits on the power to energy (P/E) ratio of of the ES (which depends on the technology adopted), on the capital available for investment in ES, and by the need to achieve a minimum rate of return $\chi$ on these investments: 
\begin{align}
\rho\up{min} e\up{R}_b \leq p\up{R}_b &\leq \rho\up{max} e\up{R}_b\,\label{PD:UL_C1}\,,\\
\textstyle\sum_{b\in B}\big(c\up{p} p\up{R}_b + c\up{e} e\up{R}_b\big) &\leq c\up{ic,max}\,,
\label{PD:UL_C2}\\
C\up{R}( x\up{P}_j, x\up{D}_j) &\geq \chi C\up{E}(x\up{U})\,,
\label{PD:UL_ESP}
\end{align}
where the ES operational profit $C\up{R}$ is calculated as follows:
\begin{gather}
\textstyle C\up{R}(x\up{P}_j, x\up{D}_j) := \sum_{j\in J}\omega_j\sumtb \big( p\up{dis}_{j,t,b}\lambda\up{lmp}_{j,t,b}\eta\up{dis}
\nonumber\\
- p\up{ch}_{j,t,b}\lambda\up{lmp}_{j,t,b} /\eta\up{ch}
+ r\up{eu}_{j,t,b} \lambda\up{eu}_{j,t}  \eta\up{dis}
+ r\up{ed}_{j,t,b} \lambda\up{ed}_{j,t} / \eta\up{ch}\nonumber\\
-c\up{dis}_{b} p\up{dis}_{j,t,b}-c\up{ch}_{b} p\up{ch}_{j,t,b}-c\up{eu}_b r\up{eu}_{j,t,b}-c\up{ed}_b r\up{ed}_{j,t,b}\big)
\,.
\end{gather}
The first two terms calculates the payment the ES receives from the energy market that settles in the locational marginal price (LMP), the third and the fourth term calculates the payment from the regulation market that settles in the system-wide regulation up and down prices, the last four terms represents the operation cost of discharging, charging, as well as providing regulation.

\subsection{Lower-Level Problem: Economic Dispatch}

Each lower-level problem minimizes the system operating cost, $C\up{P}_j$, for a particular typical day using an hourly interval. 
This economic dispatch takes into account the generation and regulation cost of conventional generators and ES units, as well as the cost associated with spillage of renewable energy. For each typical day $j$, this problem can be formulated in a compact way as follows: 
\begin{align}
&\textstyle \min_{\mathnormal{x}\up{P}_j} C\up{P}_j (\mathnormal{x}\up{P}_j) := \sumtb  c\up{rs}p\up{rs}_{j,t,b}\nonumber\\
&+ \textstyle\sumti\big(c\up{g}_i p\up{g}_{j,t,i}+c\up{gu}_i r\up{gu}_{j,t,i}+c\up{gd}_i r\up{gd}_{j,t,i}\big)\\
&+ \textstyle\sumtb\big(c\up{dis}_{b} p\up{dis}_{j,t,b}+c\up{ch}_{b} p\up{ch}_{j,t,b}+c\up{eu}_b r\up{eu}_{j,t,b}+c\up{ed}_b r\up{ed}_{j,t,b}\big)\,,\label{PD:PLL_obj}\\
&\text{subject to:} \nonumber \\
&\mathbf{M}\up{P}_j \mathnormal{x}\up{P}_j + \mathbf{M}\up{E} \mathnormal{x}\up{U} \leq \mathbf{V}\up{P}_j\,,
\label{PD:PLL_C}
\end{align}
where the decision variables are $\mathnormal{x}\up{p}_j = \{$ $p\up{ch}_{j,t,b}$,  $p\up{dis}_{j,t,b}$, $p\up{g}_{j,t,i}$, $p\up{rs}_{j,t,b}$, $r\up{ed}_{j,t,b}$, $r\up{eu}_{j,t,b}$, $r\up{gd}_{j,t,i}$, $r\up{gu}_{j,t,i}$, $e\up{soc}_{j,t,b}$, $f_{j,t,l}$, $\theta_{j,t,b}\}$. $\mathbf{M}\up{P}_j$, $\mathbf{M}\up{E}$, $\mathbf{V}\up{P}_j$ are constraint coefficient matrices. The compact expression of the constraints \eqref{PD:PLL_C} is expanded below and the dual variables associated with each constraint are shown in parentheses after a colon.
\subsubsection{Nodal power balance equations ($\forall t\in T, b\in B$)}
At each bus, the sum of the power injections and the inflows must be equal to the demand:
\begin{gather}
\textstyle \sum_{i\in I\ud{b}}p\up{g}_{j,t,i}-\sum_{l|b\in o(l)}f_{j,t,l}+\sum_{l|b\in r(l)}f_{j,t,l}
+G\up{rn}_{j,t,b}\nonumber\\-p\up{rs}_{j,t,b}+p_{j,t,b}\up{dis} \eta\up{dis} - p_{j,t,b}\up{ch}/ \eta\up{ch} = D_{j,t,b}:(\lambda\up{lmp}_{j,t,b})\label{PLL:C_bus}
\end{gather}
where $\lambda\up{lmp}_{j,t,b}$ is the locational marginal price.

\subsubsection{Regulation requirement ($t\in T, b\in B_E, i\in I$)}
Hourly up/down regulation requirements are expressed as a percentage of the system-wide demand ($\phi\up{D}$) plus a percentage of the system-wide renewable injection ($\phi\up{R}$): 
\begin{gather}
\textstyle \sum_{b\in B} r\up{eu}_{j,t,b}\eta\up{dis} + \sum_{i\in I}r\up{gu}_{j,t,i} \nonumber\\
\textstyle \geq \sum_{b\in B}\big[\phi\up{R}p\up{rn}_{j,t,b}+\phi\up{D} D_{j,t,b}\big]:(\lambda\up{ru}_{j,t})\,  \label{PLL:C_res1}\\
\textstyle \sum_{b\in B} r\up{ed}_{j,t,b}/\eta\up{ch} + \sum_{i\in I}r\up{gd}_{j,t,i} \nonumber\\
\textstyle \geq \sum_{b\in B}\big[\phi\up{R}p\up{rn}_{j,t,b}+\phi\up{D} D_{j,t,b}\big]:(\lambda\up{rd}_{j,t})\,,  \label{PLL:C_res2}
\end{gather}
where $\lambda\up{ru}_{j,t}$ and $\lambda\up{rd}_{j,t}$ are the hourly up and down regulation prices.

\subsubsection{Energy storage constraints ($\forall t\in T, b\in B_E$)} The evolution of the state of charge $e\up{soc}_{j,t,b}$ is calculated from the energy market dispatch schedules:
\begin{gather}\label{PLL_CES1}
 e\up{soc}_{j,t,b} - e\up{soc}_{j,t-1,b} = p\up{ch}_{j,t,b}-p\up{dis}_{j,t,b} :(\gamma\up{e}_{j,t,b})\,,
\end{gather}
the initial ES SoC is set to zero ($e\up{soc}_{j,0,b} = 0$) for all ES operations, and the end-of-day SoC is not enforced in \eqref{PLL_CES1}. The charging and discharging power must remain within the rated power:
\begin{gather}
p\up{ch}_{j,t,b} + r\up{ed}_{j,t,b}\leq p\up{R}_b:(\varphi\up{ch}_{j,t,b})\,\label{PLL_CES2} \\
p\up{dis}_{j,t,b} + r\up{eu}_{j,t,b}\leq p\up{R}_b:(\varphi\up{dis}_{j,t,b})\,,\label{PLL_CES3} 
\end{gather}
and the ES must sustain the full regulation reserve dispatch for the required time interval ($T\up{es}$):
\begin{gather}
e\up{soc}_{j,t,b} +T\up{es} r\up{ed}_{j,t,b}\leq e\up{R}_b:(\varphi\up{soc}_{j,t,b})\,\label{PLL_CES5}\\
e\up{soc}_{j,t,b}- T\up{es} r\up{eu}_{j,t,b}\geq 0:(\psi\up{soc}_{j,t,b})\,.\label{PLL_CES4}
\end{gather}

\subsubsection{Other constraints} Appendix~\ref{App_SL} defines the formulation of the generator power rating and ramp constraints, the constrains on renewable spillage, and the dc power flow model used to enforce the network constraints.

\subsection{Dual Lower-Level Problem}
We apply the primal-dual transformation to the primal lower-level~(PLL) problem due to its convexity. The dual lower-level~(DLL) problem optimizes system prices so that constraint \eqref{PD:UL_ESP} is enforced.
The DLL objective function is formulated as follows: 
\begin{align}
&\textstyle \max_{\mathnormal{x}\up{D}_j} C\up{D}_j(\mathnormal{x}\up{U}, \mathnormal{x}\up{D}_j) := \textstyle
\sumti\big[(\varphi\up{g}_{j,t,i} G\up{max}_{i}+\psi\up{g}_{j,t,i} G\up{min}_{i})\nonumber\\
&\textstyle
+ R\up{u}_i (\varphi\up{R}_{j,t,i}+T\up{ru} \varphi\up{gu}_{j,t,i})+R\up{d}_i (-\psi\up{R}_{j,t,i}+T\up{rd} \varphi\up{gd}_{j,t,i})\big]\nonumber\\
&\textstyle
+\sum_{i\in I}(\varphi\up{R}_{e,1,i}+\varphi\up{R}_{e,1,i}) G_{j,t}\up{0}+\sum_{\subalign{t&\in T \\ l&\in L}}(\varphi\up{f}_{j,t,l}-\psi\up{f}_{j,t,l}) F\up{max}_l\nonumber\\
&\textstyle
+ \sumtb \Big[ \varphi\up{rs}_{j,t,b}G\up{rs}_{j,t,b}+ \lambda_{j,t,b}\up{lmp}(D_{j,t,b}-G\up{rn}_{j,t,b})\big]\nonumber\\
&\textstyle
+\sumtb \big[\phi\up{rn} G\up{rn}_{j,t,b}(\lambda_{j,t}\up{ru}+\lambda_{j,t}\up{rd}) +\phi\up{D} D_{j,t,b}(\lambda_{j,t}\up{ru}+\lambda_{j,t}\up{rd})\big]\nonumber\\&\textstyle
+\sumtb \big[p\up{R}_b(\varphi\up{ch}_{j,t,b}+\varphi\up{ dis}_{j,t,b})+e\up{R}_b\varphi\up{soc}_{j,t,b}\big]\,.\\
&\text{subject to:} \nonumber \\
&\mathbf{M}\up{D}_j\mathnormal{x}\up{D}_j \leq \mathbf{V}\up{D}_j\,.
\end{align}
where $\mathbf{M}\up{D}_j$ and $\mathbf{V}\up{D}_j$ are the constraint coefficient matrices. The detail of these constraints are given in Appendix~\ref{App_SL}.
Note that the objective function includes products of UL ($x\up{U}_j$) and DLL variables ($x\up{D}_j$).

\begin{figure}[t]
	\centering
	\includegraphics[trim = 55mm 75mm 43mm 0mm, clip, width = .9  \columnwidth]{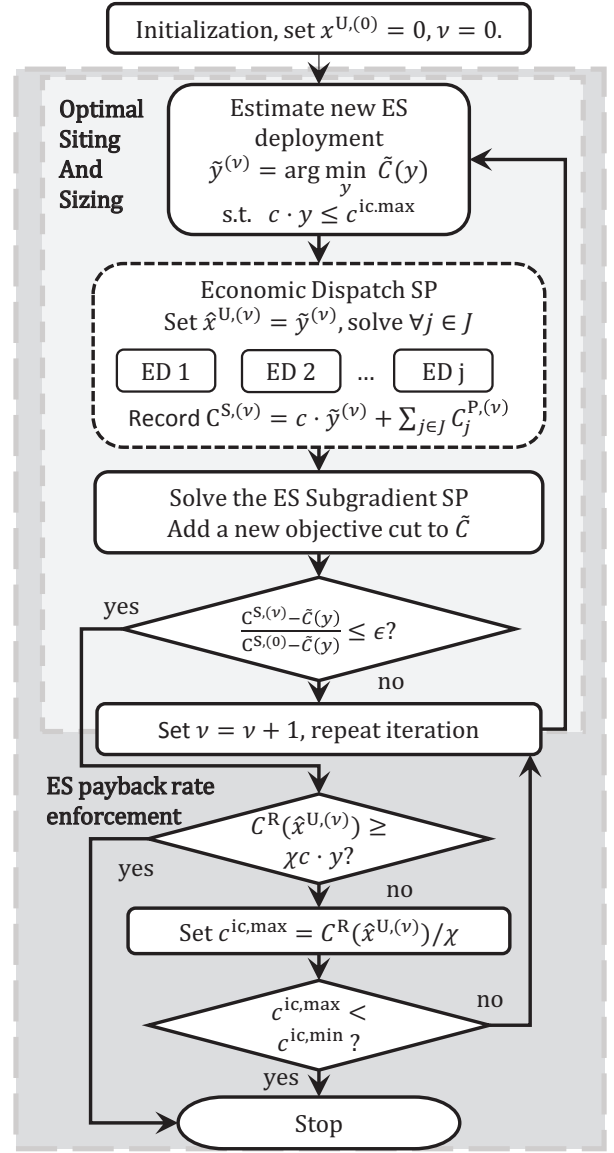}
	\caption{Flowchart of the solution algorithm 
	}
	\label{fig:cp}
\end{figure}
 
\section{Solution Method}\label{Sec:Meth}

Decomposition has been used extensively for solving large-scale programming problems~\cite{zhao2014marginal,linderoth2003decomposition,birge1985decomposition,nielsen1997scalable}, especially for scenario-based stochastic programmings~\cite{nap2016next, papavasiliou2011reserve, kazempour2012strategic,baringo2012wind,nasrolahpour2016strategic}. A stochastic planning problem couples independent scenarios with a few planning decision variables. An effective decomposition breaks each scenario into a subproblem, which can be solved in sequence or in parallel. It is also easier to aggregate subproblem results in such decomposition structures, and the master problem can be solved rapidly and accurately.
Fig.~\ref{fig:cp} illustrates the proposed solution algorithm that involves an inner-loop and an outer-loop. The inner-loop identifies the optimal ES locations subject to the maximum ES investment budget, $c\up{ic,max}$, and the outer-loop enforces the ES rate of return constraint \eqref{PD:UL_ESP}.

In the inner-loop, the main problem is decomposed into scenario subproblems by fixing the value of ES planning variables. Each scenario subproblem solves an ED problem. The inner-loop is initialized with no ES installation in the system, and solve the ED for all scenarios. Based on the ED results, the potential benefit of ES installation is calculated for each bus in the system in a subgradient form. ES installations are updated accordingly. Therefore, the inner-loop calculates ES siting and sizing decisions iteratively, until the estimated distance to the exact optimal solution is sufficiently small. The decomposition technique is described in Section~\ref{Sec:PDc}. Section~\ref{Sec:ESSG} explains how the subgradient of the objective function with respect to ES planning variables is calculated. Section~\ref{Sec:CP} explains the subgradient cutting-plane method used to update ES planning variables at each iteration and solve the optimal ES location problem.

In the outer-loop, the optimal ES siting and sizing decisions are tested against the ES rate of return constraint \eqref{PD:UL_ESP}. If it is not satisfied, the maximum ES investment budget, $c\up{ic,max}$, is reduced and the inner-loop is repeated (see Section~\ref{Sec:ESP}). The algorithm terminates once a current solution satisfy constraint \eqref{PD:UL_ESP} or the maximum ES investment budget, $c\up{ic,max}$, reaches the minimum ES investment limit, $c\up{ic,min}$.

\subsection{Problem Decomposition}\label{Sec:PDc}
The bi-level problem \eqref{PD:UL_Obj}--\eqref{PLL_CES4} can be recast into a single-level (SL) equivalent. The  objective function of this problem is:
\begin{align}
&\hspace{-0.2cm} \min_{\mathnormal{x}\up{U}, \mathnormal{x}\up{P}_j, \mathnormal{x}\up{D}_j} C\up{S}(\mathnormal{x}\up{U}, \mathnormal{x}\up{P}_j) :=\sum_{j\in J}\omega_jC\up{P}_j(\mathnormal{x}\up{U}, \mathnormal{x}\up{P}_j) + C\up{E}(\mathnormal{x}\up{U})\,,
\label{SL:obj}\\
&\text{subject to:} \nonumber\\
&\text{UL, PLL, and DLL constraints} \\
&C\up{P}_j(\mathnormal{x}\up{U}, \mathnormal{x}\up{P}_j) = C\up{D}_j(\mathnormal{x}\up{U}, \mathnormal{x}\up{D}_j)\,,\;j\in J\,,
\label{SL:SDC}
\end{align}
where \eqref{SL:SDC} represents the strong duality constraint. The details of the formulation of this problem are given in Appendix~\ref{App_SL}. 

When the value of the coupling variables $\mathnormal{x}\up{U}$ is fixed, we can then apply primal decomposition to this problem. For the sake of simplicity, we first ignore the profit constraint \eqref{PD:UL_ESP}. The subproblem becomes a linear ED problem for each typical day. This decomposed SL problem can then be solved iteratively as follows~\cite{boyd2004convex}:
\subsubsection{Set initial values for the coupling variables} The solution algorithm starts with $\mathnormal{x}\up{U,(0)}=0$, indicating no ES deployment. 
\subsubsection{Solve the subproblems} At iteration $\nu$, set $\mathnormal{x}\up{U}=\mathnormal{x}\up{U,(\nu)}$, solve each EDSP in parallel to obtain $\hat{\mathnormal{x}}\up{P, (\nu)}_j$ and $\hat{\mathnormal{x}}\up{D, (\nu)}_j$. 
\subsubsection{Solve the master problem} Calculate the subgradients of $C\up{S}$ with respect to $x\up{U}$ and update the UL variables accordingly.
\subsubsection{Iteration} Check for convergence, and repeat from Step 2) if needed. The convergence criterion is explained in Section.~\ref{Sec:CP}.

While subgradient methods can be easily used to solve the master problem, their convergence is slow and they do not provide a measurement of the optimality of the results. We therefore incorporate the subgradient cutting-plane method in the proposed approach because it converges faster and has a non-heuristic stopping criterion.

\subsection{The Cutting-plane Method}\label{Sec:CP}

We apply cutting-plane methods~\cite{kelley1960cutting, belloni2005introduction,  boyd2007localization} to solve the master problem. Cutting-plane methods incorporate results from previous iterations and form a piece-wise linear approximation ($\tilde{C}^{(\nu)}$) of the objective function:
\begin{gather}
\textstyle \tilde{C}^{(\nu+1)}(y) := \max_{k\leq \nu} [C\up{S}(\tilde{y}\up{(\nu)})+(y-\tilde{y}\up{(\nu)})\cdot g\up{U, (\emph{k})}]\,,
\label{CP_obj}
\end{gather}
where $y$ is an inquiry point identical to $x\up{U}$, and $\tilde{C}^{(\nu)}(\tilde{y}\up{(\nu)})$ is a lower-bound estimate of the optimal objective function value, i.e., the system cost.

At each iteration, the ES subproblems are solved by setting $x\up{U,(\nu)}=\tilde{y}\up{(\nu)}$ where $\tilde{y}\up{(\nu)}\in \arg \min_y \tilde{C}^{(\nu)}(y)$ subject to constraints \eqref{PD:UL_C1} and \eqref{PD:UL_C2}. The inquiry point $\tilde{y}\up{(\nu)}$, the current system cost value, and the calculated subgradients are then added to \eqref{CP_obj} as a new objective cut for future iterations. As iterations proceed, $\tilde{C}^{(\nu)}$ approaches the actual objective function, and a solution is found when the difference between $C\up{S}(x\up{U,(\nu)})$ and $\tilde{C}^{(\nu)}(\tilde{y}\up{(\nu)})$ is below a  tolerance. We assume that the algorithm terminates based on a relative tolerance with respect to the estimated maximum system cost saving: 
\begin{gather}
\tilde{C}^{(\nu)}(\tilde{y}\up{(\nu)})-C\up{S}(x\up{U,(\nu)}) \leq \epsilon[C\up{S}(0)-\tilde{C}^{(\nu)}(\tilde{y}\up{(\nu)})]\,,
\end{gather}
where $C\up{S}(0)$ is the system cost without ES deployments, and $\epsilon$ is the relative tolerance. Therefore the system cost saving in the algorithm is always greater than $(1-\epsilon)$ of the optimal system cost saving. 

\subsection{Energy Storage Subgradient Cuts  }\label{Sec:ESSG}
ES subgradients for $B\up{E}$ buses are calculated directly using dual variables associated with constraints \eqref{PLL_CES2}--\eqref{PLL_CES4}, this derivation is shown in Appendix~\ref{App:SBE}. However, the value of $\hat{\varphi}\up{soc, (\nu)}_{j,t,b}$ is always zero at the first iteration, because energy rating constraints are implicit, thus \eqref{PLL_CES5} never binds when the ES power rating is zero. Moreover, $\hat{\varphi}\up{ch, (\nu)}_{j,t,b}$ and $\hat{\varphi}\up{dis, (\nu)}_{j,t,b}$ cannot reflect the value for arbitrage unless the ES has non-zero SoC evolutions to correlate temporal arbitrage decision between different time intervals. Therefore, we calculate ES subgradients for $B\up{N}$ buses through a subgradient subproblem (SGSP) to correctly identify the value of marginal ES investment. The SGSP provides an initial ES P/E ratio $\rho_b^0$, so that solutions to the master problem has reduced perturbations and converges faster. The derivation of the SGSP is shown in Appendix~\ref{App:SBN}.
At iteration $\nu$, the subgradients of $C\up{S}$ with respect to the ES power rating ($g\up{p}_b$) and energy rating ($g\up{e}_b$) are calculated as:
\begin{align} 
g\up{p,(\nu)}_b &=
\begin{cases}
\;\, c\up{p}+\hat{\varphi}\up{ch, (\nu)}_{j,t,b} + \hat{\varphi}\up{dis, (\nu)}_{j,t,b}\,, & \quad b\in B\up{E}\,\\
\;\, \hat{g}\up{0,(\nu)}_b \hat{\rho}^{0,(\nu)}_b/(1+\hat{\rho}^{0,(\nu)}_b)  \,, & \quad b\in B\up{N}\,
\end{cases}\label{ESSG:P}
\\
g\up{e,(\nu)}_b &= 
\begin{cases}
\;\, c\up{e}+\hat{\varphi}\up{soc, (\nu)}_{j,t,b}\,,  &  \qquad\quad b\in B\up{E}\,\\
\;\, \hat{g}\up{0,(\nu)}_b /(1+\hat{\rho}^{0,(\nu)}_b)  \,, &    \qquad\quad b\in B\up{N}\,,
\end{cases}\label{ESSG:E}
\end{align}
where $\hat{g}\up{0, (\nu)}_b$ and $\hat{\rho}^{0, (\nu)}_b$ are determined by solving the following subgradient subproblem (SGSP):
\begin{align}
\textstyle &\textstyle \min_{p\up{ch}_{j,t,b}, p\up{dis}_{j,t,b}, r\up{eu}_{j,t,b}, r\up{ed}_{j,t,b}, e\up{soc}_{j,t,b}, \rho^0_b} g\up{0,(\nu)}_b :=  \nonumber\\
\textstyle &\textstyle\quad\sum_{j\in J}\omega_j\sum_{t\in T}\big[ p\up{ch}_{j,t,b}\hat{\lambda}\up{lmp, (\nu)}_{j,t,b}/\eta\up{ch}
-p\up{dis}_{j,t,b}\hat{\lambda}\up{lmp, (\nu)}_{j,t,b}\eta\up{dis} \nonumber\\
&\textstyle\quad- r\up{eu}_{j,t,b}\hat{\lambda}\up{ru, (\nu)}_{j,t}\eta\up{dis} - r\up{ed}_{j,t,b}\hat{\lambda}\up{rd, (\nu)}_{j,t}/\eta\up{ch} + c\up{dis}_{b} p\up{dis}_{j,t,b}\nonumber\\ 
&\textstyle\quad +c\up{ch}_{b} p\up{ch}_{j,t,b}+c\up{eu}_b r\up{eu}_{j,t,b}+c\up{ed}_b r\up{ed}_{j,t,b}\big] + \rho\up{0}_bc\up{p}+c\up{e}\,,
\label{ESSG:BN}
\end{align}
subject to constraints \eqref{PD:UL_C1} and \eqref{PLL_CES1}--\eqref{PLL_CES4} by setting $p\up{R}_b = \rho^0_b$ and $e\up{R}_b=1$. This subproblem maximizes the profit of incremental ES deployments at $B\up{N}$ buses, where ES are price-takes and profit maximization is equivalent to system operating cost minimization~\cite{castillo2013profit}. $\hat{\rho}^0_b$ is the optimal P/E ratio for price-taker ES deployments, which is close to the true optimal P/E ratio if the ES has limited price influences. Subgradients at $B\up{N}$ buses are designed to enforce $\hat{\rho}^0_b$ over all new ES deployments. ES deployments with near-optimal P/E ratios have faster convergence due to minimum perturbations between ES power and energy investment decisions.

\subsection{Incorporating the ES Profit Constraint}\label{Sec:ESP}

In Appendix~\ref{App:ESP}, we show that the ES operational revenue can be represented using the ES subgradients:
\begin{align}
C\up{R} = -\textstyle\sum_{b\in B}\big[(g\up{p}_b-c\up{p})p\up{R}_b + (g\up{e}_b-c\up{e})e\up{R}_b\big]\,.
\end{align}
Because all ES allocation variables must have non-negative values, all $p\up{R}_b$ and $e\up{R}_b$ with non-zero values must have negative subgradients. we can therefore infer that $C\up{R} \geq C\up{E}$ in all optimal locations, and in unconstrained ES locations, the ES rate of return $\chi$ converges to one. Hence $\chi \geq 1$ is guaranteed for all optimal ES locations. For $\chi > 1$, we can reasonably assume that ES has a limited effect on system prices and that the system-wide ES operating revenue should only increase with ES investment. Therefore, for the optimal ES locations, the ES revenue $\hat{C}\up{R}$ is a concave monotonic increasing function of the ES investment cost $C\up{E}$ such that: 
\begin{gather}
0 \leq d\hat{C}\up{R}(C\up{E})/dC\up{E} \leq 1\,,
\label{ESPA_A}
\end{gather}
where  $C\up{E}$ is capped by $c\up{ic,max}$ in constraint \eqref{PD:UL_C2}. If a rate of return $\chi$ is achievable in the system, then there must be some $C\up{E'}$ that satisfy:
\begin{gather}
\hat{C}\up{R}(C\up{E'}) - \chi C\up{E'} \geq 0\,.
\end{gather}
When an ES investment cost $C\up{E}$ violates \eqref{PD:UL_C2}, we can estimate an upper-bound of $C\up{E'}$ as
\begin{gather}
C\up{E'} \leq \hat{C}\up{R}(C\up{E'})/\chi \leq \hat{C}\up{R}(C\up{E})/\chi\,,\label{ESP:UB}
\end{gather}
because $C\up{R}(C\up{E'}) \leq C\up{R}(C\up{E})$ according to \eqref{ESPA_A}. Therefore, \eqref{PD:UL_ESP} can be satisfied by iteratively solving the optimal ES allocation with a reduced maximum ES investment cost $c\up{ic,max}=C\up{R}/\chi$. 

Since we use the cutting-plane method to solve the master problem and the feasible region is reduced when setting $c\up{ic,max}=C\up{R}/\chi$ (recall that \eqref{PD:UL_ESP} only binds when $\chi>1$), solving the optimal ES allocations recursively will not add much complexity because the method already has a fairly good estimate of the objective function.

\subsection{Comparison to Benders Decomposition}

Benders decomposition is a classic approach for solving block-structured optimization with coupling (complicating) variables~\cite{conejo2006decomposition}, and has been extensively used for solving strategic bi-level planning problems in power system~\cite{kazempour2012strategic,baringo2012wind,nasrolahpour2016strategic}. The proposed algorithm is similar to Benders decomposition because it decomposes the optimization problem by fixing the coupling variables and solves the master problem using cutting planes. However, it incorporates two key improvements over a classic Benders decomposition.

First, it uses coordinated subgradient cuts, which provides more accurate information on the value of marginal ES investments than Benders dual cuts. Since the value of an ES deployment is jointly affected by its power and energy ratings, these planning decisions are not independent, and the P/E ratio must be optimized. However, Benders dual cuts based on the binding conditions of ES rating constraints \eqref{PLL_CES2}~-~\eqref{PLL_CES4} are inefficient at coordinating investments on power and energy ratings, especially during the early iterations where ES ratings are mostly zero. These uncoordinated cuts cause the master problem solution to oscillate around the optimal point, and significantly slow down the convergence. Instead of using dual cuts, the proposed algorithm uses coordinated ES subgradients, which enforces near-optimal P/E ratios over all new ES deployments and thus speeds up the convergence.

Besides using coordinated subgradient cuts, we analytically derived a relationship between the maximum ES investment budget and the profitability, and decomposed the bi-level problem into a recursive structure, as shown in Fig.~\ref{fig:cp}. Previous studies~\cite{kazempour2012strategic,baringo2012wind,nasrolahpour2016strategic} follow an approach that combines Benders decomposition with MPEC. This approach first materializes the bi-level problem into a MPEC problem, the MPEC problem is then recast into a MILP structure using the 'big M' method~\cite{floudas1995nonlinear}, and Benders decomposition is applied to decompose and solve the MILP problem. The 'big M' method uses auxiliary integer variables and a sufficiently large constant $M$ to linearize nonlinear terms. However the accuracy and computational speed of the 'big M' linearization are very sensitive to the value of $M$, if $M$ is not large enough, the linearization is not accurate, if $M$ is too large, the computation can be extremely slow. Compared to the MPEC+Benders method, our algorithm requires no auxiliary linearization variables and the size of the master problem does not increase with the number of scenarios. Therefore, the proposed algorithm has better scalability and leads to more robust planning result. In addition, the proposed algorithm generates simpler subproblems that are solved faster than in the MPEC+Benders method. In Section~\ref{Sec:CS:ESP:CTC}, we  demonstrate that the computational speed of the proposed method surpasses the solution time achieved in a similar previous study.

\section{Case Study Test System}\label{Sec:Case}

\subsection{System Settings}

The proposed ES planning model and solution method were tested using a modified 240-bus reduced WECC system~\cite{price2011reduced}. This system includes 448 transmission lines, 71 aggregated thermal plants and renewable sources including hydro, wind, and solar. 
The maximum expected forecasts of all renewable generations are grouped as $G\up{rn}_{j,t,b}$, the maximum allowable spillage for hydro generation is enforced in $G\up{rs}_{j,t,b}$, other types of renewable generation have no curtailment limits. Renewable curtailments are necessary in the modified WECC testbed because large renewable generation capacities are installed at some buses with a limited transmission capacity, and the objective of the economic dispatch is to minimize the system operating cost. In certain cases, such as days with strong winds and low demand, a certain amount of wind power generation must be curtailed to maintain secure operation of the system.

We use a `3+5'\% reserve policy for setting the requirements for regulation~\cite{papavasiliou2011reserve}, hence  $\phi\up{D}=3\%$ and $\phi\up{R}=5\%$. Regulation parameters are adjusted so that the regulation prices are identical to the actual day-ahead regulation clearing prices in  CAISO~\cite{CAISO_Report2014}. 
The value of renewable spillages is set to zero. The modified WECC system has a daily ED operating cost ranging from 15 to 35 M\$.

All simulations were carried out in CPLEX under GAMS~\cite{GAMS} on an Intel Xenon 2.55 GHz processor with 32 GB of RAM. Typical days and their respective weights from the year-long demand and renewable generation profiles were identified using a hierarchical clustering algorithm~\cite{pitt2000applications}. The convergence criterion is set to $\epsilon=5\%$.

\begin{table}[t]
	\begin{center}
		\centering
		\caption{Energy Storage Model Parameters.}
		\label{tab:ES_tech}
		\begin{tabular}{r || c || c}
			\hline
			\hline
			Energy storage technology & AA-CAES & LiBES 	\Tstrut\Bstrut\\
			\hline
			\hline	
			Power rating investment (\$/kW-year) & 1250 - 20 & 409 - 20 	\Tstrut\Bstrut\\
			\hline
			Energy rating investment (\$/kWh-year) & 150 - 20 & 468 - 20	\Tstrut\Bstrut\\
			\hline
			Battery replacement cost (\$/kWh) & n/a & 406 	\Tstrut\Bstrut\\
			\hline
			P/E ratio range (h$^{-1}$) & 0.05 to 0.25 & 0.1 to 4 	\Tstrut\Bstrut\\
			\hline
			Round trip efficiency & 0.72 & 0.9 	\Tstrut\Bstrut\\
			\hline	
			Incremental production cost (\$/MWh) & 0 & 87 \Tstrut\Bstrut\\
			\hline
			Incremental consumption cost (\$/MWh) & 0 & 0 \Tstrut\Bstrut\\
			\hline
			\hline						
		\end{tabular}
	\end{center}
	\vspace{-.5cm}
\end{table}

\subsection{Energy Storage Cost Model}\label{Sec:Ben:ESC}

We consider two types of representative ES technologies: 1) above-ground advanced adiabatic compressed air energy storage (AA-CAES), and 2) lithium-ion battery energy storage (LiBES)~\cite{zakeri2015electrical, nykvist2015rapidly}. The ES cost model consists of three parts: 
\begin{itemize}
    \item Investment cost of the power equipments ($c\up{p}$) proportional to the ES power rating (unit: \$/kW). This investment covers the turbine generator and air compressor in AA-CAES, or the power electronic equipments in LiBES.
    \item Investment cost of the storage system ($c\up{e}$) proportional to the ES energy rating (unit: \$/kWh). This investment covers the storage tank in AA-CAES, or the battery management system in LiBES.
    \item Marginal production cost (unit: \$/kWh). This is the energy production cost of ES units. 
\end{itemize}
Because we consider large-scale ES installations, we assume that the fixed storage installation cost, including the land and construction cost, scales linearly with the ES rating and can therefore be incorporated in the costs proportional to the power and energy ratings.
AA-CAES units have high investment cost for power ratings and low investment cost for the storage capacity, the operation cost is also negligible because AA-CAES consumes no fuel for power generation. Electrochemical battery energy storage such as LiBES~\cite{dunn2011electrical, zakeri2015electrical} has more evenly distributed investment cost components. The lifetime of lithium batteries is very sensitive to operations due to degradation, the marginal production cost of LiBES is therefore defined based on its cell cycle life.

\begin{figure}
    \centering
    \includegraphics[trim = 5mm 0mm 10mm 0mm, clip, width = 1\columnwidth]{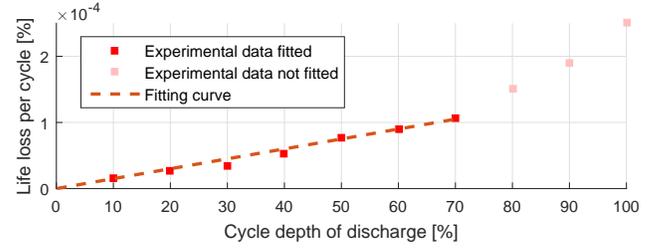}
    \caption{MO battery cycle life curve and fitting to 70\% DoD.}
    \vspace{0mm}
    \label{Fig:ES_deg}
\end{figure}

We assume that cycle aging only occurs during battery discharging and has a constant marginal cost. Using cycle life test data set for Lithium manganese oxide (LMO) batteries~\cite{xu2016modeling}, we apply a linear fit to the LMO cycle life loss per cycle up to 70\% depth of discharge (DoD) (Fig.~\ref{Fig:ES_deg}). The operation region of the LiBES is limited to the range from 20\% to 90\% of SoC to avoid deep discharges as well as overcharge and overdischarge effects, because these factors severely reduce the battery life~\cite{vetter2005ageing}. Instead of introducing new SoC constraints, the LiBES is oversized to reflect the increased cost due to a narrowed SoC operation region. We assume that the battery cells in the LiBES are always replaced once reaching their end of life. The marginal discharging cost of LiBES is calculated by prorating the battery cell replacement cost to the cycle life loss curve
\begin{align}
    c\up{dis} = a\up{fit}\frac{\text{Cell replacement cost (\$406/kWh)}}{\text{DoD operation range (70\%)}}\,,
\end{align}
where $a\up{fit}$ is the linear fitted slop of the cycle life loss curve in Fig.~\ref{Fig:ES_deg}.

Table~\ref{tab:ES_tech} shows the capital cost, P/E ratio range, efficiency, and operating cost of these two ES cost models. The battery cells in LiBES are replaced once reaches the end of life. We assume a 5\% annual interest rate and calculate the daily prorated cost as in~\cite{pandzic2015near}. All 240 buses are considered as ES deployment candidates. 

\subsection{Negative Pricing and Storage Dispatch}

In optimal ED solutions, an ES unit may be dispatched to charge and discharge simultaneously during the occurrence of negative LMPs~\cite{go2016assessing}. The storage round-trip efficiency causes energy spillages that are beneficial to the system, and ES units gain additional revenue. Such dispatches are physically achievable for AA-CAES units because the air compressor and generator use separate pipelines~\cite{steta2010modeling}, so that the compressor and generator can operate at the same time. 

LiBES can only charge or discharge at one time and simultaneous dispatches must be avoided. In Appendix~\ref{App:SCER} we demonstrated the following sufficient condition for avoiding simultaneous charging and discharging ($\forall b\in B, t\in T, j\in J$):
\begin{align}
    c\up{dis} + c\up{ch} > -({1}/{\eta\up{dis}}-\eta\up{ch})\lambda\up{lmp}_{j,t,b}\,. 
\end{align}
The above sufficient condition explains that in optimal ED solutions, an ES unit will not charge and discharge simultaneously as long as the operation cost for a round-trip dispatch is higher than the product of the round trip efficiency loss and the negative LMP value. In other words, the cost of performing simultaneous dispatches is higher than the market payment. In the modified WECC model, the largest negative LMP never exceeds -200~\$/MWh, and the round-trip efficiency of LiBES is 90\%. Therefore as long as the marginal production cost of LiBES is higher than 20~\$/MWh, simultaneous LiBES dispatches are avoided.

\begin{figure}%
	\centering
	\subfloat[Daily prorated system cost saving.]{
		\includegraphics[trim = 109mm 31mm 108mm 25mm, clip, width = .5\columnwidth]{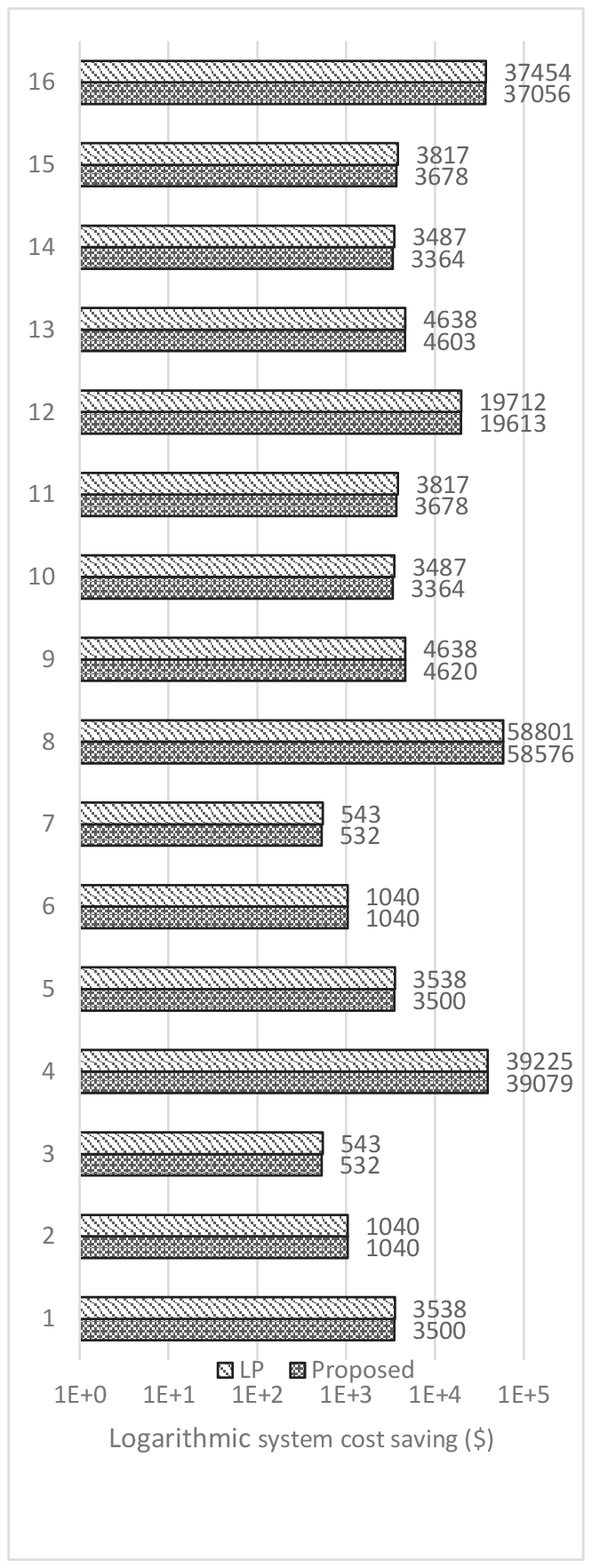}
		\label{Fig:Comp_SC}%
	}
	\subfloat[Computation time.]{
		\includegraphics[trim = 109mm 31mm 108mm 25mm, clip, width = .5\columnwidth]{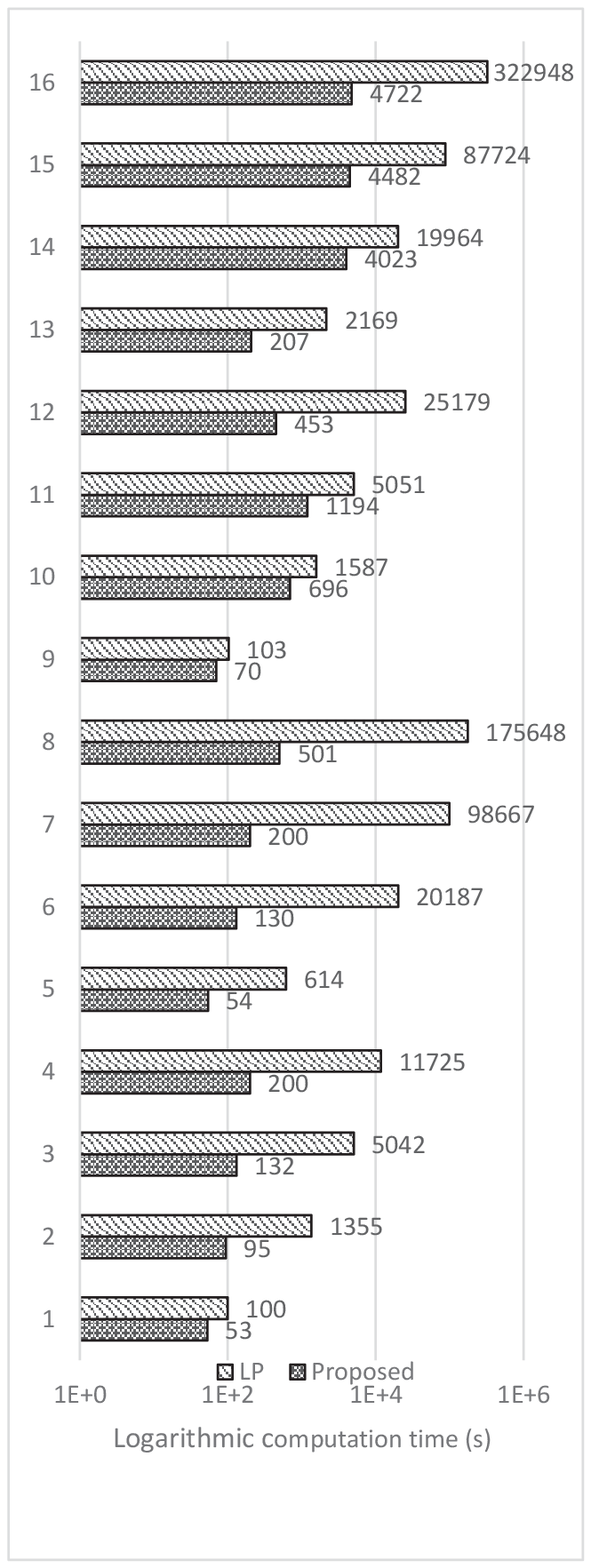}
		\label{Fig:Comp_TC}%
	}
	\caption{Computation test results of the 16 test cases.}%
	\label{Fig:Comp}%

\end{figure}

\subsection{Regulation Cost and Dispatch Model}

The cost of providing regulation is estimated using normalized CAISO area control error (ACE) data~\cite{caiso_ace, xu2016comparison}. The provision of regulation does not increase the operating cost for AA-CAES because these units have no marginal operating cost. The average cost for LiBES to provide 1 MW of regulation up for one hour is 10\% of the marginal production cost ($c\up{eu} = 0.1c\up{dis}$) under the assumption that an average 100 kWh of energy is generated. The regulation cost scales linearly with the regulation up capacity within the 70\% DoD region. The LiBES has no cost for providing regulation down because charging has no marginal cost.

In CAISO, ES units have two options for participating regulation: the regulation energy management program (REM) and the traditional option (non-REM). In the REM program, CAISO co-optimize the ACE with the real-time energy market and generates a regularized regulation signal that has a 15-minute zero-mean in energy. Therefore REM units are only required to have a 15 minute capacity ($T\up{es}=0.25$) and are not deviated from their scheduled SoC levels.

Non-REM units are required to have a continuous full dispatch time requirement of one hour ($T\up{es}=1$). SoC deviations for providing regulation in Non-REM ES units are not accounted in the ED formulation because we primarily evaluate the economic value of ES investment in hourly economic dispatch. In the regulation market, energy deviations caused by regulation provision are settled at the real-time locational marginal prices after the dispatch period~\cite{caiso_settlement}. Therefore, from an economic point of view, ES does not gain or loses energy in regulation provision (i.e., ES cannot receive ‘free energy for charging’ by providing regulation, the charged energy is still settled at market prices), and the proposed planning model leads to sufficiently accurate decisions without considering the real-time regulation energy deviations. In real-time dispatches, ES units can adopt control strategies against large energy deviations~\cite{xu2014bess} and maintain the scheduled dispatch.

\section{Simulation Results}\label{Sec:Test}

\subsection{Computational Performance}

\begin{table*}[t]
	\centering
	\caption{Comparison of ES rate of return with a maximum ES investment cost of 50~k\$/day.}
	\label{tab:ES_pc}
	\begin{tabular}{r || c c c c|| c c c c|| c c c c }
	\hline
	\hline
	Selected typical day & \multicolumn{4}{c||}{Day 100} & \multicolumn{4}{c||}{Day 141} & \multicolumn{4}{c}{Day 285} \Tstrut\\
	\hline
	\hline
	ES rate of return (\%) & 100 & 110 & 120 & 150 & 100 & 110 & 120 & 150 & 100 & 110 & 120 & 150 \Tstrut\Bstrut\\
	\hline
	Runtime (s) & 102 & 164 & 148 & 219 & 86 & 91 & 91 & 113 & 398 & 639 & 551 & 566 \Tstrut\Bstrut\\
	\hline
	ED operation cost (M\$/day) & 15.83 & 15.84 & 15.85 & 15.86 & 22.87 & 22.87 & 22.87 & 22.91 & 25.19 & 25.20 & 25.20 & 25.20\Tstrut\Bstrut\\
	\hline
	ES operation revenue (k\$/day) & 27.2 & 13.6 & 10.0 & 0 & 60.0 & 60.0 & 60.0 & 36.4 & 9.2 & 1.9 & 0 & 0 \Tstrut\Bstrut\\
	\hline
	ES investment cost (k\$/day) & 26.8 & 12.3 & 8.0 & 0 & 50.0 & 50.0 & 50.0 & 24.1 & 9.1 & 1.7 & 0 & 0 \Tstrut\Bstrut\\
	\hline 
	ES location (bus number) & 155 & 155 & 155 & n/a & 155 & 155 & 155 & 155 & 228 & 228 & n/a & n/a \Tstrut\Bstrut\\
	\hline
	\hline
	\end{tabular}
\end{table*}

We compare the computational performance of the proposed method against solving the problem directly using CPLEX. When the profit constraint \eqref{PD:UL_ESP} is ignored, the objective function \eqref{SL:obj} and constraints \eqref{PD:UL_C1}, \eqref{PD:UL_C2} and \eqref{PD:PLL_C} become a linear problem (LP), and can be solved by using the solver CPLEX. 

We designed 16 test cases with different planning scenarios. Case 1-4 are the optimal ES allocation considering 1, 3, 5, and 10 typical days, subject to a maximum ES investment budget constraint, using  the CAES ES model. Case 5-8 are identical to 1-4 except that they do not include a maximum investment constraint. Case 9-16 are identical to 1-8 except that the LiBES ES model is used.

As shown in Figure.~\ref{Fig:Comp}, the proposed method is significantly faster than solving the LP problem directly using CPLEX, while all system cost saving results are within the set tolerance. CPLEX exhibits an approximately quadratic increase in computation time with the number of typical days, while the proposed method demonstrates a much slower increase. However, the computational speed of the proposed method depends on the renewable and demand profiles because some profiles result in smoother system cost functions, which facilitates the convergence of the subgradient cutting-plane method. 
When the IC constraint is excluded, the search region for ES allocation expands and thus the computation time of both methods increases. However, this effect is much smaller in the proposed method.

\subsection{Rate of Return on ES Investment}
We performed ES planning on three different days subject to different ES rate of return constraints. Table~\ref{tab:ES_pc} shows the results for the LiBES model with $T\up{es}=0.25$. A higher ES rate of return reduces the installed ES capacity and increases the system operating cost. A return rate of 150\% is only achievable in one of the three days.

The computation time of the proposed method increases moderately when the payback rate is greater than 1, because the optimal ES allocation is solved repeatedly. However this will not result in a polynomial or exponential increase in complexity because the cutting-plane method keeps track of historical results. Enforcing a higher rate of return reduces the maximum ES investment budget and hence decreases the feasible region. In turn this reduces the solution time when the problem is solved iteratively. 

This table also shows that buses 155 and 285 are the only locations where ES is deployed for day 100, 141 and 285. In other single-day tests that we performed, ES was also located at buses 15, 90, 198, 226, 227, 228. These buses are good locations for performing spatio-arbitrage because they are connected to frequently congested lines and renewable sources, especially hydro units. In particular, LMPs are frequently negative at bus 155.  

\begin{figure}[t]%
	\centering
	\subfloat[Installed power rating.]{
		\includegraphics[trim = 5mm 0mm 10mm 0mm, clip, width = .95\columnwidth]{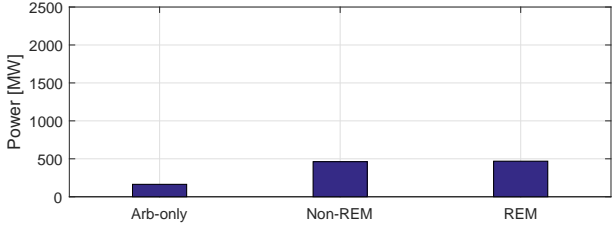}
		\label{Fig:CAES_P}%
	}
	\\
	\subfloat[Installed energy rating.]{
		\includegraphics[trim = 5mm 0mm 10mm 0mm, clip, width = .95\columnwidth]{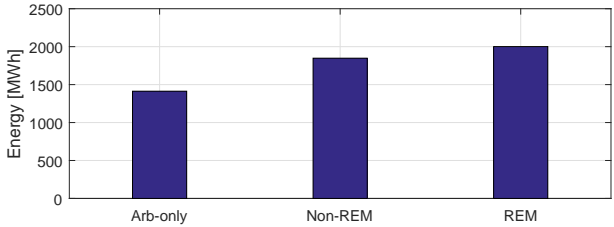}
		\label{Fig:CAES_E}%
	}
	\\
	\subfloat[Annual market revenue.]{
		\includegraphics[trim = 5mm 0mm 10mm 0mm, clip, width = .95\columnwidth]{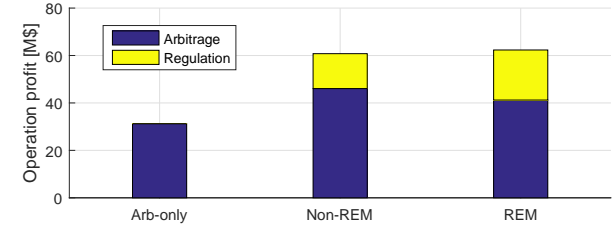}
		\label{Fig:CAES_rev}%
	}
	\caption{AA-CAES planning in different market scenarios.}%
	\label{Fig:CAES}
\end{figure}

\begin{figure}[t]%
	\centering
	\subfloat[Installed power rating.]{
		\includegraphics[trim = 5mm 0mm 10mm 0mm, clip, width = .95\columnwidth]{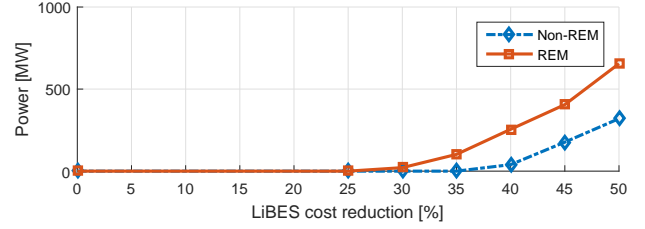}
		\label{Fig:LiBES_P}%
	}
	\\
	\subfloat[Installed energy rating.]{
		\includegraphics[trim = 5mm 0mm 10mm 0mm, clip, width = .95\columnwidth]{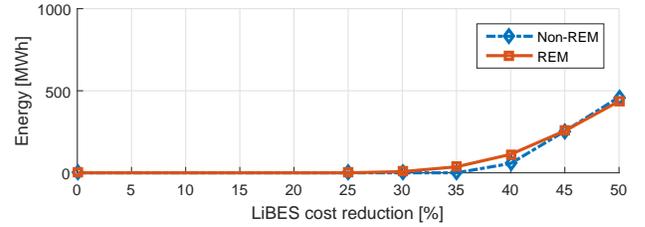}
		\label{Fig:LiBES_E}%
	}
	\\
	\subfloat[Market revenue from arbitrage and regulation.]{
		\includegraphics[trim = 5mm 0mm 10mm 0mm, clip, width = .95\columnwidth]{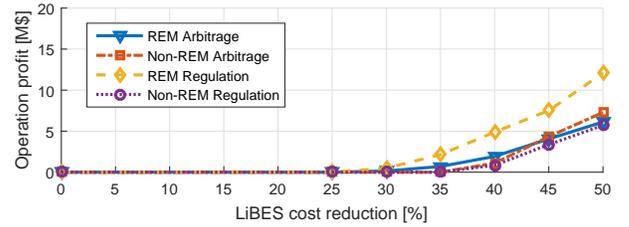}
		\label{Fig:LiBES_rev}%
	}
	\caption{LiBES planning with decreasing investment cost.}%
	\label{Fig:LiBES}
\end{figure}

\subsection{Stochastic ES Planning}

We performed stochastic ES planning considering 20 typical days with no maximum ES investment limit. Three market scenarios are considered: in \emph{Arb-only} AA-CAES only participates in the energy market, in \emph{Non-REM} AA-CAES can participate energy and regulation markets under traditional regulation requirements, and in \emph{REM} AA-CAES can participate energy and regulation markets under REM regulation requirements. 

\subsubsection{AA-CAES Results} Fig.~\ref{Fig:CAES} shows stochastic planning results for AA-CAES.
Because system-wide regulation prices are independent of the location, arbitrage is the sole factor for ES siting. Bus 155 is the optimal choice for all ES allocations, mainly due to its high occurrence of negative LMPs. In the Arb-only market scenario, AA-CAES has a P-to-E ratio of 0.11, equivalent to 9 hours of rated energy capacity. In Non-REM and REM cases, the planning results have larger installation capacities, and the P-to-E ratio also increases. The change in the regulation requirement does not have a significant impact on the planning results, and arbitrage is still the primary market income source.

\subsubsection{LiBES Results}
No LiBES is installed in any market scenarios under current investment cost as shown in Table~\ref{tab:ES_tech}. Since the decreasing trend of LiBES investment cost is expected to continue for the next ten years~\cite{dunn2011electrical,nykvist2015rapidly},  it is reasonable to assess the planning of LiBES using reduced investment cost. In Fig.~\ref{Fig:LiBES}, LiBES planning results are shown for up to 50\% investment cost reduction for the Non-REM and REM market scenarios, while in the Arb-only case, no LiBES is installed. The result shows that investments in LiBES will become profitable when the investment cost dropped by at least 30\% from its current value , and the installed capacity increases steadily with further cost reductions. At each cost level, the market revenue from arbitrage is roughly the same in the Non-REM and the REM case, while the regulation revenue almost doubles in REM. The difference in revenue also reflects in the installation capacity, while the installed energy capacity is similar in the REM and Non-REM cases, LiBES has a much higher P-to-E ratio with REM.

\subsection{Computational Speed Comparison}\label{Sec:CS:ESP:CTC}

The proposed solution algorithm is faster than solving a bi-level ES planning problem using the combination of Benders decomposition and MPEC.
Nasrolahpour~\emph{et al.}~\cite{nasrolahpour2016strategic} use the MPEC+Benders approach to solve a bi-level ES planning problem for optimal ES sizing in a single bus system in an energy market environment. Their solution time range from 3 to 6 hours. By comparison, when our method was applied to the optimization of ES siting and sizing in a 240-bus system considering both energy and regulation markets, the longest simulation finished within one hour, and most simulations finished within 15 minutes.

\section{Summary}\label{Sec:Con}

In this paper, we have formulated the optimal ES profit-constrained siting and sizing as a bi-level problem with a minimum rate of return constraint. We have proposed a scalable solution method involving a primal decomposition and subgradient cutting-planes. The proposed method is significantly faster than CPLEX for solving LP ES planning problems.

The proposed solution method has the same order of complexity as conventional economic dispatch, thus making this method computationally tractable for any system with a feasible ED solution. Since the decomposed subproblems are independent of each other, the computation time increases linearly with the number of typical days considered. The solution time could be further improved by solving the subproblems in parallel.

We have analyzed the optimal ES siting in joint energy and reserve markets on a modified WECC 240-bus model. The sensitivity of these siting decisions has been studied with respect to different ES technologies, the rate of return on ES investments, and regulation market policies. The results show that increasing the rate of return requirement greatly reduces the deployment of ES.
In the stochastic ES planning, AA-CAES shows a higher potential for reducing system cost than LiBES, which depends  on the design of the regulation market for its profitability. However AA-CAES technology is still at the pilot stage, while grid-scale installations of LiBES are happening worldwide.

\section*{Acknowledgement}
The  authors  thank  Dr.  I.  Gyuk  and  his  colleagues  at  the US  DOE  Energy  Storage  Program  for  funding  this  research. Sandia  National  Laboratories  is  a  multi-program  laboratory managed and operated by Sandia Corporation, a wholly owned subsidiary  of  Lockheed  Martin  Corporation,  for  the  U.S.Department  of  Energy’s National  Nuclear  Security  Administration under Contract DE-AC04-94AL85000. 

\appendix

\subsection{Single-Level Equivalent Problem Formulation}\label{App_SL}

This problem consists of the objective function \eqref{SL:obj}, the UL constraints \eqref{PD:UL_C1}--\eqref{PD:UL_ESP}, and the following constraints:
\subsubsection{PLL constraints} Equations \eqref{PLL:C_bus}--\eqref{PLL_CES4} and the following constraints ($\forall i\in I, j\in J, t\in T, b\in B$):
\begin{gather}
G\up{min}_i + r\up{gd}_{j,t,i} \leq p\up{g}_{j,t,i}\leq G\up{max}_i-r\up{gu}_{j,t,i}:(\psi\up{g}_{j,t,i},\varphi\up{g}_{j,t,i})\,\label{SL_C1b}\\
0 \leq r\up{gu}_{j,t,i}\leq T\up{ru} R\up{u}_i:(\psi\up{gu}_{j,t,b}, \varphi\up{gu}_{j,t,b})\,\label{SL_C1c}\\
0 \leq r\up{gd}_{j,t,i}\leq T\up{rd} R\up{d}_i:(\psi\up{gd}_{j,t,b}, \varphi\up{gd}_{j,t,b})\,\label{SL_C1d}\\
-R\up{d}_i \leq p\up{g}_{j,t,i}-p\up{g}_{j,t-1,i}\leq R\up{u}_i:(\psi\up{R}_{j,t,i}, \varphi\up{R}_{j,t,i})\,\label{SL_C1e}\\
f_{j,t,l}=(\theta_{j,t,o(l)}-\theta_{j,t,r(l)})/x_l : (\gamma\up{f}_{j,t,l})\,\label{SL_C2a}\\
-F\up{max}_l\leq f_{j,t,l} \leq F\up{max}_l : (\psi\up{f}_{j,t,l}, \varphi\up{f}_{j,t,l})\,\label{SL_C2b}\\
0\leq p\up{rs}_{j,t,b} \leq G\up{rs}_{j,t,b}:(\psi\up{rn}_{j,t,b}, \varphi\up{rn}_{j,t,b})\,,\label{SL_C3}
\end{gather}
The minimum and maximum capacity of all generators are enforced in \eqref{SL_C1b}. \eqref{SL_C1c} and \eqref{SL_C1d} model the ramp requirement for regulation, and \eqref{SL_C1e} models the ramp requirement for dispatch. The DC power flow is modeled in \eqref{SL_C2a} and \eqref{SL_C2b}. The maximum expected forecast for renewable generation and the maximum allowable renewable spillage are enforced in \eqref{SL_C3}.
\subsubsection{DLL constraints} 
The DLL problem has the following constraints ($\forall i\in I, j\in J, t\in T, b\in B$):
\begin{align}
\varphi_{j,t,l}\up{f}+\psi\up{f}_{j,t,l}+\gamma\up{f}_{j,t,l}-\lambda_{j,t,o(l)}\up{lmp}+\lambda_{j,t,r(l)}\up{lmp} &=0\,\label{SL_DC1}\\
\psi\up{rn}_{j,t,b}+\varphi\up{rn}_{j,t,b}-\lambda_{j,t,b}\up{lmp}+(\lambda_{j,t}\up{ru} + \lambda_{j,t}\up{rd})\phi\up{D} &=c\up{rs}\,\label{SL_DC2}\\
\varphi\up{g}_{j,t,i}+\psi\up{g}_{j,t,i}+\varphi\up{R}_{j,t,i}-\varphi\up{R}_{j,t+1,i}+\psi\up{R}_{j,t,i}-&\psi\up{R}_{j,t+1,i}\nonumber\\
+\lambda_{j,t,b(i)}\up{lmp}&=c_i\up{g}\,\label{SL_DC3}\\
\varphi\up{g}_{j,t,i}+\lambda_{j,t}\up{ru}+\varphi\up{gu}_{j,t,i}+\psi\up{gu}_{j,t,i}&=c\up{gu}_i\,\label{SL_DC5}\\
-\psi\up{g}_{j,t,i} +\lambda_{j,t}\up{rd}+\varphi\up{gd}_{j,t,i}+\psi\up{gd}_{j,t,i}&=c\up{gd}_i\,\label{SL_DC6}\\
\varphi\up{soc}_{j,t,b}+\psi\up{soc}_{j,t,b}+\gamma\up{soc}_{j,t,b}-\gamma\up{soc}_{j,t+1,b}&=0\,\label{SL_DC7}\\
\varphi\up{soc}_{j,n_T,b}+\psi\up{soc}_{j,n_T,b}+\gamma\up{soc}_{e,n_T,b} &= 0\,\label{SL_DC8}\\
\varphi_{j,t,b}\up{ch}+\psi_{j,t,b}\up{ch}-\gamma\up{soc}_{j,t,b}-\lambda\up{lmp}_{j,t,b}/\eta\up{ch} &= c\up{ch}\,\label{SL_DC9}\\
\varphi_{j,t,b}\up{dis}+\psi_{j,t,b}\up{dis}+\gamma\up{soc}_{j,t,b}+ \lambda\up{lmp}_{j,t,b} \eta\up{dis} &= c\up{dis}\,\label{SL_DC10}\\
\varphi_{j,t,b}\up{ch}+\psi_{j,t,b}\up{rd}+T\up{es}\varphi_{j,t,b}\up{soc} - \overline{e}\up{ed}\gamma\up{soc}_{j,t,b}+ \lambda\up{rd}_{j,t}/\eta\up{ch} &= c\up{ed}\,\label{SL_DC11}\\
\varphi_{j,t,b}\up{dis}+\psi_{j,t,b}\up{ru}-T\up{es}\psi_{j,t,b}\up{soc} - \overline{e}\up{eu}\gamma\up{soc}_{j,t,b} +\lambda\up{ru}_{j,t}\eta\up{dis}  &= c\up{eu}\,,\label{SL_DC12}
\end{align}
where $\psi \geq 0$ and $\varphi \leq 0$.

\subsection{ES Profit Constraint Transformation}\label{App:ESP}
From DLL constraints \eqref{SL_DC9}--\eqref{SL_DC10} and \eqref{SL_DC11}--\eqref{SL_DC12}, we obtain the following equalities:
\begin{align}
&p\up{dis}_{j,t,b}(\lambda\up{lmp}_{j,t,b}\eta\up{dis} - c\up{dis})
- p\up{ch}_{j,t,b}(\lambda\up{lmp}_{j,t,b} /\eta\up{ch} + c\up{ch})
\nonumber\\
&=\gamma\up{soc}_{j,t,b}(p\up{ch}_{j,t,b}-p\up{dis}_{j,t,b}) 
- p\up{dis}_{j,t,b}(\psi\up{dis}_{j,t,b}
+ \varphi\up{dis}_{j,t,b})
\nonumber\\&-p\up{ch}_{j,t,b}(\psi\up{ch}_{j,t,b} 
+ \varphi\up{ch}_{j,t,b})\,\label{ESP_Eq2}\\
&r\up{eu}_{j,t,b} (\lambda\up{eu}_{j,t}\eta\up{dis} - c\up{eu})
+ r\up{ed}_{j,t,b}(\lambda\up{ed}_{j,t} / \eta\up{ch} - c\up{ed}) \nonumber\\ &= 
-(\varphi_{j,t,b}\up{ch}+\psi_{j,t,b}\up{rd} +T\up{es}\varphi_{j,t,b}\up{soc}) r\up{ed}_{j,t,b} 
\nonumber\\&-(\varphi_{j,t,b}\up{dis}+\psi_{j,t,b}\up{ru}-T\up{es}\psi_{j,t,b}\up{soc}) r\up{eu}_{j,t,b}\,.
\label{ESP_Eq3}
\end{align}

By using \eqref{PLL_CES1}, \eqref{SL_DC7}, and \eqref{SL_DC8}, we derive the following expression:
\begin{align}
&\textstyle \sum_{t\in T}\gamma\up{soc}_{j,t,b}(p\up{ch}_{j,t,b}-p\up{dis}_{j,t,b})  \nonumber\\
&=\textstyle \sum_{t=1}^{n_{T}-1}e\up{soc}_{j,t,b}(\gamma\up{soc}_{j,t,b} - \gamma\up{soc}_{j,t+1,b}) + e\up{soc}_{j, n_T,b}\gamma\up{soc}_{j, n_T, b}\nonumber\\
&\textstyle= \sum_{t\in T}e\up{soc}_{j,t,b}(\psi\up{soc}_{j,t,b}+\varphi\up{soc}_{j,t,b})\,,
\label{ESP_Eq4}
\end{align}
We can obtain the linear daily revenue collected by ES by 1) combining and rearranging \eqref{ESP_Eq2}--\eqref{ESP_Eq4}, and 2) substituting the complementary slackness condition associated with \eqref{PLL_CES2}--\eqref{PLL_CES4}. 
\begin{align}
&(p\up{ch}_{j,t,b} + r\up{ed}_{j,t,b})\varphi\up{ch}_{j,t,b} + (p\up{dis}_{j,t,b} + r\up{eu}_{j,t,b})\varphi\up{dis}_{j,t,b} \nonumber\\
&+ (e\up{soc}_{j,t,b}+T\up{es}r\up{ed}_{j,t,b})\varphi\up{soc}_{j,t,b} + (e\up{soc}_{j,t,b}-T\up{es}r\up{eu}_{j,t,b})\psi\up{soc}_{j,t,b} \nonumber\\
&=  \textstyle\sum_{t\in T} \big[(\varphi\up{ch}_{j,t,b} + \varphi\up{dis}_{j,t,b}) p\up{R}_b + \varphi\up{soc}_{j,t,b} e\up{R}_b\big]\,,
\label{ESP_Eq5}
\end{align}
which leads to
\begin{align}
C\up{R} = -\textstyle\sum_{j\in J}\sumtb\big[(\varphi\up{ch}_{j,t,b} + \varphi\up{dis}_{j,t,b}) p\up{R}_b + \varphi\up{soc}_{j,t,b} e\up{R}_b\big]\,.
\label{ESP_Eq6}
\end{align}
Because ES revenue only applies to $b\in B\up{E}$, by comparing \eqref{ESP_Eq6} with \eqref{ESSG:P} and \eqref{ESSG:E},we can represent ES revenue using ES subgradients:
\begin{align}
C\up{R} = -\textstyle\sum_{b\in B}\big[(g\up{p}_b-c\up{p})p\up{R}_b + (g\up{e}_b-c\up{e})e\up{R}_b\big]\,.
\end{align}

\subsection{ES Subgradient Derivation at $B\up{E}$ Buses}\label{App:SBE}
We calculate the ES subgradients assuming $C\up{S}$ is differentiable. At $(\nu)$th iteration, the ES subgradients $g\up{U, (\nu)}$ includes the subgradients with respect to $p\up{R}_b$ and $e\up{R}_b$ for $b\in B$
\begin{gather}
g\up{U,(\nu)} = [g\up{p,(\nu)}_b\quad  g\up{e,(\nu)}_b]^T\,.
\end{gather}
$g\up{U,(\nu)}$ can be calculated using either the primal or the dual form of the ED problem with their minimizer (or maximizer):
\begin{align}
g\up{U, (\nu)} \approx &\textstyle \nabla_{\mathnormal{x}\up{U}}C\up{S}(\mathnormal{x}\up{U, (\nu)}, \hat{\mathnormal{x}}\up{P, (\nu)}) \\
=&\nabla_{\mathnormal{x}\up{U}}C\up{E}(\mathnormal{x}\up{U, (\nu)})
+\textstyle \sum_{j\in J}\omega_j\nabla_{\mathnormal{x}\up{U}}C\up{P}_j(\mathnormal{x}\up{U, (\nu)}, \hat{\mathnormal{x}}\up{P, (\nu)}_j)\nonumber\\
=&\nabla_{\mathnormal{x}\up{U}}C\up{E}(\mathnormal{x}\up{U, (\nu)})
+\textstyle \sum_{j\in J}\omega_j\nabla_{\mathnormal{x}\up{U}}C\up{D}_j(\mathnormal{x}\up{U, (\nu)}, \hat{\mathnormal{x}}\up{D, (\nu)}_j)\,.\nonumber
\end{align}
and the subgradients are calculated as follows:
\begin{align}
\textstyle \lim_{\Delta x\up{U} \to 0} &||C\up{D}_j(\mathnormal{x}\up{U, (\nu)}+\Delta x\up{U}, \hat{\mathnormal{x}}\up{D, (\nu)}_j) 
- C\up{D}_j(\mathnormal{x}\up{U, (\nu)}, \hat{\mathnormal{x}}\up{D, (\nu)}_j)\nonumber\\ 
&- [g\up{U}]^T\Delta x\up{U} || / ||\Delta x\up{U}||= 0\,.
\end{align}
We use the dual form of the ED problem, and the subgradient for $b\in B\up{E}$ is:
\begin{align}
\begin{bmatrix}
g\up{p, (\nu)} \\[1em] g\up{e, (\nu)}
\end{bmatrix}
=
\begin{bmatrix}
c\up{p}+\sum_{j\in J}\omega_j\sumtb(\hat{\varphi}\up{ch, (\nu)}_{j,t,b} + \hat{\varphi}\up{dis, (\nu)}_{j,t,b})\\[1em]
 c\up{e}+\sum_{j\in J}\omega_j\sumtb\hat{\varphi}\up{soc, (\nu)}_{j,t,b}
\end{bmatrix}\,.
\end{align}

\subsection{ES Subgradient Derivation for $B\up{N}$ Buses}\label{App:SBN}
For $b\in B\up{N}$, let $\Delta p\up{R}_b$ and $\Delta e\up{R}_b$ be sufficiently small. We use the strong duality condition and replace $C\up{P}_j$ with $C\up{D}_j$ in $C\up{S}$. Because $\Delta p\up{R}_b\to 0$, $\Delta e\up{R}_b\to 0$, other decision variables are not affected and are removed, leaving only terms that are directly associated with energy storage at $B\up{N}$ buses and obtain the following problem that calculates the ES gradient at $B\up{N}$ buses at  iteration $\nu$:
\begin{gather}
\textstyle \max_{x\up{\Delta}} C\up{0, (\nu)}(\hat{x}\up{\Delta}, \hat{x}\up{D}_j) := \nonumber\\
\textstyle \sum_{j\in J}\omega_j\sumtb \big[\Delta p\up{R}_b({\varphi}\up{ch}_{j,t,b}+{\varphi}\up{ dis}_{j,t,b})+\Delta e\up{R}_b{\varphi}\up{soc}_{j,t,b}\big]
\end{gather}
we let $\rho^0_b = \Delta p\up{R}_b / \Delta e\up{R}_b$, because all $\Delta p\up{R}_b$ and $\Delta e\up{R}_b$ variables are completely independent, thus the problem is equivalent to:
\begin{align}
&\textstyle \max_{\Delta x} \sum_{j\in J}\omega_j\sumtb \big[\rho^0_b({\varphi}\up{ch}_{j,t,b}+{\varphi}\up{ dis}_{j,t,b})+{\varphi}\up{soc}_{j,t,b}\big] \\
&\text{subject to:} \nonumber \\
&\rho\up{min} \leq \rho^0_b \leq \rho\up{max}\,,\label{g0_C1}
\end{align}
and \eqref{SL_DC5} to \eqref{SL_DC12} by replacing $\lambda\up{lmp}_{j,t,b}$, $\lambda\up{ru}_{j,t,b}$, $\lambda\up{rd}_{j,t,b}$ with $\hat{\lambda}\up{lmp, (\nu)}_{j,t,b}$, $\hat{\lambda}\up{ru, (\nu)}_{j,t,b}$, $\hat{\lambda}\up{rd, (\nu)}_{j,t,b}$ in $x\up{D, (\nu)}_j$ because these ES do not affect prices.
This subproblem can be transformed into its equivalent primal form
\begin{gather}
\textstyle \min_{x\up{\Delta}} C\up{0}(x\up{\Delta}, \hat{x}\up{D, (\nu)}_j) := \sum_{b\in B} g\up{0}(x\up{\Delta}, \hat{x}\up{D, (\nu)}_j)\,,
\end{gather}
which is equivalent to \eqref{ESSG:BN}. 

\subsection{Exact Relaxation of ES Dispatch Constraints}\label{App:SCER}

In the established ED problem, an ES can be enforced to only charge or discharge at a single time step with the following non-convex complementary constraint~\cite{li2016sufficient} ($\forall$ $j\in J$, $t\in T$, $b\in B$)
\begin{align}
    p\up{dis}_{j,t,b}p\up{ch}_{j,t,b} = 0\,,\label{SCER:cc}
\end{align}
Sufficient conditions for an exact relaxation of \eqref{SCER:cc} is analyzed using the Karush-Kuhn-Tucker~(KKT) condition. 
In the KKT condition for the ED problem in \eqref{PD:PLL_obj} and \eqref{PD:PLL_C}, the derivative of the Lagrangian function with respect to ES discharging variables $p\up{dis}_{j,t,b}$ must equal to zero, hence the following equation holds ($\forall$ $j\in J$, $t\in T$, $b\in B$)
\begin{align}
    c\up{dis} - \varphi\up{dis}_{j,t,b} - \psi\up{dis}_{j,t,b} + \gamma\up{e}_{j,t,b} + \lambda\up{lmp}_{j,t,b}/\eta\up{dis} = 0
    \label{SCER:dis}
\end{align}
similarly, for ES charging variables $p\up{ch}_{j,t,b}$, the following equation holds ($\forall$ $j\in J$, $t\in T$, $b\in B$)
\begin{align}
    c\up{ch} - \varphi\up{ch}_{j,t,b} - \psi\up{ch}_{j,t,b} - \gamma\up{e}_{j,t,b} - \lambda\up{lmp}_{j,t,b}\eta\up{ch} = 0
    \label{SCER:ch}
\end{align}
Assume there exists $p\up{dis}_{j,t,b} > 0 $ and $p\up{ch}_{j,t,b} > 0$ at bus $b$ at time $t$ during typical day $j$ in the optimal solution of the ED problem. Then $\psi\up{dis}_{j,t,b}=0$, $\psi\up{ch}_{j,t,b}=0$ because of the complementary slackness conditions. Summing \eqref{SCER:dis} and \eqref{SCER:ch} and the following equation holds
\begin{align}
    c\up{dis} + c\up{ch} -\varphi\up{dis}_{j,t,b}-\varphi\up{ch}_{j,t,b} + ({1}/{\eta\up{dis}}-\eta\up{ch})\lambda\up{lmp}_{j,t,b} = 0\,, \label{SCER:e1}
\end{align}
because $\varphi\up{dis}_{j,t,b} \leq 0$, $\varphi\up{ch}_{j,t,b} \leq 0$, \eqref{SCER:e1} can be reduced to
\begin{align}
    c\up{dis} + c\up{ch} \leq -({1}/{\eta\up{dis}}-\eta\up{ch})\lambda\up{lmp}_{j,t,b}\,, \label{SCER:e2}
\end{align}
\eqref{SCER:e2} describes the necessary condition for $p\up{dis}_{j,t,b} > 0 $ and $p\up{ch}_{j,t,b} > 0$. Hence, the sufficient condition for the exact relaxation of the complementary constraint of \eqref{SCER:cc} is 
\begin{align}
    c\up{dis} + c\up{ch} > -({1}/{\eta\up{dis}}-\eta\up{ch})\lambda\up{lmp}_{j,t,b}\,, \label{SCER:e3}
\end{align}
for all $j\in J$, $t\in T$, $b\in B$.

\bibliographystyle{IEEEtran}	
\bibliography{IEEEabrv,literature}		

\begin{IEEEbiographynophoto}{Bolun Xu}
(S'14) received B.S. degrees in Electrical and Computer Engineering
from the University of Michigan, Ann Arbor, USA and Shanghai Jiaotong
University, Shanghai, China in 2011, and the M.Sc degree in Electrical
Engineering from Swiss Federal Institute of Technology, Zurich, Switzerland
in 2014.

He is currently pursuing the Ph.D. degree in Electrical Engineering at the
University of Washington, Seattle, WA, USA. His research interests include
energy storage, power system operations, and power system economics.
\end{IEEEbiographynophoto}

\begin{IEEEbiographynophoto}{Yishen Wang}
(S'12) received the B.S. degree from the Department of
Electrical Engineering, Tsinghua University, Beijing, China, in 2011. He is
currently pursuing the Ph.D. degree in electrical engineering at the University
of Washington, Seattle, WA, USA.

His research interests include power system economics and operation,
energy storage, renewable forecasting and electricity markets.
\end{IEEEbiographynophoto}

\begin{IEEEbiographynophoto}{Yury Dvorkin}
(S'11-M'16) received his Ph.D. degree from the University of
Washington, Seattle, WA, USA, in 2016.

Dvorkin is currently an Assistant Professor in the Department of Electrical
and Computer Engineering at New York University, New York, NY, USA.
Dvorkin was awarded the 2016 Scientific Achievement Award by Clean
Energy Institute (University of Washington) for his doctoral dissertation “Operations
and Planning in Sustainable Power Systems”. His research interests
include short- and long-term planning in power systems with renewable
generation and power system economics.
\end{IEEEbiographynophoto}

\begin{IEEEbiographynophoto}{Ricardo~Fern\'andez-Blanco}
(S'10-M'15) received the Ingeniero Industrial
degree and the Ph.D. degree in electrical engineering from the Universidad
de Castilla-La Mancha, Ciudad Real, Spain, in 2009 and 2014, respectively.

He is currently a Scientific/Technical Project Officer at the JRC.C7 Knowledge
for the Energy Union (Joint Research Center), Petten, The Netherlands.
His research interests include the fields of operations and economics of power
systems, bilevel programming, hydrothermal coordination, and electricity
markets.
\end{IEEEbiographynophoto}

\begin{IEEEbiographynophoto}{Jean-Paul Watson}
(M'10) received the B.S., M.S., and Ph.D. degrees in
computer science.

He is a Distinguished Member of Technical Staff with the Discrete Math
and Optimization Department, Sandia National Laboratories, Albuquerque,
NM, USA. He leads a number of research efforts related to stochastic
optimization, ranging from fundamental algorithm research and development,
to applications including power grid operations and planning.
\end{IEEEbiographynophoto}

\begin{IEEEbiographynophoto}{Cesar A. Silva-Monroy}
(M'15) received the B.S. degree from the Universidad Industrial de Santander, Bucaramanga, Colombia, and the M.S. and Ph.D. degrees from University of Washington, Seattle, WA, USA, all in electrical engineering. 

He is a Senior Member of Technical Staff with the Electric Power Systems Research Department at Sandia National Laboratories, Albuquerque, NM, USA. He leads and collaborates in several projects that seek to increase the efficiency, resilience, and reliability of the grid by applying advanced optimization techniques to power system operations, planning, and control.
\end{IEEEbiographynophoto}

\begin{IEEEbiographynophoto}{Daniel S. Kirschen}
(M’86-SM’91-F’07) received his electrical and mechanical engineering degree from the Universite Libre de Bruxelles, Brussels, Belgium, in 1970 and his M.S. and Ph.D. degrees from the University of Wisconsin, Madison, WI, USA, in 1980, and 1985, respectively.

He is currently the Donald W. and Ruth Mary Close Professor of Electrical Engineering at the University of Washington, Seattle, WA, USA. His research interests include smart grids, the integration of renewable energy sources in the grid, power system economics, and power system security.
\end{IEEEbiographynophoto}

\end{document}